\documentstyle[amstex,amssymb,11pt]{article}

  \input xypic
  \xyoption{curve}
  \input cyracc.def

  \newfam\cyifam
  \font\tencyi=wncyi10
  \font\sevencyi=wncyi7
  \font\fivecyi=wncyi5
  
  \textfont\cyifam=\tencyi \scriptfont\cyifam=\sevencyi
    \scriptscriptfont\cyifam=\fivecyi

  \newfam\cyrfam
  \font\tencyr=wncyr10
  \font\sevencyr=wncyr7
  \font\fivecyr=wncyr5
  
  \textfont\cyrfam=\tencyr \scriptfont\cyrfam=\sevencyr
    \scriptscriptfont\cyrfam=\fivecyr
 \newcommand{\lon}{\longrightarrow}
 
 \newcommand{\rar}{\rightarrow}
 \newcommand{\hook}{\hookrightarrow}
 \newcommand{\Proof}{{\bf Proof}.\, }
 \newcommand{\OM}{\Omega M}
 \newcommand{\LM}{\Lambda_M}
 \newcommand{\CP}{{\Bbb C} {\Bbb P}}

 \newcommand{\bS}{{\Bbb S}}
 \newcommand{\p}{{\partial}}

 \newcommand{\R}{{\Bbb R}}
 \newcommand{\ot}{\otimes}

 \newcommand{\Id}{\mbox{Id}}
 \newcommand{\Beq}{\begin{equation}}
 \newcommand{\Eeq}{\end{equation}}
 \newcommand{\Beqr}{\begin{eqnarray}}
 \newcommand{\Eeqr}{\end{eqnarray}}
 \newcommand{\Beqrn}{\begin{eqnarray*}}
 \newcommand{\Eeqrn}{\end{eqnarray*}}
 \newcommand{\Ba}{\begin{array}}
 \newcommand{\Ea}{\end{array}}
 \newcommand{\Bi}{\begin{itemize}}
 \newcommand{\Ei}{\end{itemize}}
 \newcommand{\Bc}{\begin{center}}
 \newcommand{\Ec}{\end{center}}

 \newcommand{\f}{{\cal O}}

 \newcommand{\caD}{{\cal D}}
 
 \newcommand{\cF}{{\cal F}}

 \newcommand{\cP}{{\cal P}}

 \newcommand{\al}{\alpha}
 
 \newcommand{\ga}{\gamma}
 
 \newcommand{\Ga}{\Gamma}

 \newcommand{\om}{\omega}

 \newcommand{\Ker}{{\mathsf Ker}\, }
 \newcommand{\Img}{{\mathsf Im}\, }
 
 \newcommand{\Hom}{{\mathrm Hom}}
 \newcommand{\bH}{{\bf H}}

 \newcommand{\sip}{\smallskip}
 \newcommand{\bip}{\bigskip}

\begin{document}

 \sloppy

 \title{De Rham model for string topology}
 \author{ S.A.\ Merkulov}
 \date{}
\maketitle


\bip

\bip

{\bf 0. Introduction.} Let $M$ be a simply connected closed
oriented $n$-dimensional manifold, and $LM:={\mathrm
Map}_{C^\infty}(\bS^1, M)$ the associated free loop space. String
topology of Chas and Sullivan \cite{CS} deals with  an ample
family of algebraic operations  on the ordinary and equivariant
homologies of $LM$, the most important being a graded commutative
associative product,
$$
\Cap: \bH_p(LM)\ot \bH_q(LM) \lon \bH_{p+q}(LM),
$$
on the shifted ordinary homology\footnote{In this paper all
homologies are assumed to be over $\R$.},
$\bH_\bullet(LM):=H_{\bullet + n}(LM)$. The product $\Cap$ is
compatible with (in general, non-commutative)  Pontrjagin product,
$$
\star: H_p(\OM) \ot H_q(\OM) \lon H_{p+q}(\OM),
$$
in the homology of the subspace, $\OM$, of $LM$ consisting of
loops based at some fixed point in $M$, and the same product
$\Cap$ is also
 an extension
of the classical intersection product in $H_\bullet(M)$: the
natural chain of linear maps
$$
 \left(H_\bullet(M),\cap \right)
 \stackrel{\mathrm  constant\atop loops}{\lon}
 \left(\bH_\bullet(LM),\Cap \right)
 \stackrel{\mathrm  intersection\ with \atop the\ subspace\ \OM}{\lon}
\left(H_\bullet(\OM),\star \right)
 $$
is  a chain of  homomorphisms of algebras.

\bip

Let $(\LM=\bigoplus_{i=0}^n\LM^i, d)$ be the De Rham differential
graded (dg, for short) algebra of $M$. With the help of simplicial
techniques it was established in \cite{J} that
$$
\hspace{4cm} H_\bullet(LM)= {\mathrm
Ext}_{\LM\ot\LM^{o\!p}}^\bullet\left(\LM,\LM^*\right),
\hspace{4cm} (*)
$$
i.e.\ that morphisms in the derived category, $\caD(\LM)$, of
$\LM$-bimodules between the two objects, $\LM$ and its dual
$\LM^*:=\Hom_\R(\LM,\R)$, have a nice free loop space
interpretation. As $M$ is closed oriented, the Poincare duality
asserts that $\LM^*$ and $\LM[n]$ are isomorphic in $\caD(\LM)$ so
that $(*)$ can be rewritten equivalently as
$$
 \bH_\bullet(LM) ={\mathrm Ext}_{\LM\ot\LM^{o\!p}}^\bullet\left(\LM,\LM\right).
$$
Moreover, as it was  shown in \cite{CJ,Co} using ring spectra and
simplicial methods, the latter isomorphism is actually an
isomorphism of algebras,
$$
\hspace{2cm}\left( \bH_\bullet(LM),  \Cap\right) =\left({\mathrm
Ext}_{\LM\ot\LM^{o\!p}}^\bullet\left(\LM,\LM\right), {\mathrm
Yoneda\ product}\right), \hspace{2cm} (**)
$$
i.e.\ the Chas-Sullivan product is the same thing as the
composition in the endomorphism ring of the object $\LM$ in the
derived category of $\LM$-bimodules. The isomorphism $(**)$ was
also studied
in the thesis of Tradler \cite{Tr}.

\bip

This paper
 offers  new geometrically transparent (and down-to-earth) proofs
 of the string topology main theorems
which are based on the theory of iterated integrals \cite{C}. In particular,
the isomorphism $(*)$ gets incarnated as the holonomy map. Which
combined with the Thom class of $M$  explains  why $(**)$ is an
isomorphism of algebras. As a by-product, the resulting  De Rham
model for string topology  provides us with new  algorithms for
computing both the free loop space homology and the Chas-Sullivan
product;
 see also \cite{CJY, FTV,R} for other approaches.

\bip

In fact, our approach to string topology easily generalizes to its
``brane'' version. Let $f: Z\rar M$  be a smooth map from a
compact oriented $p$-dimensional manifold to a compact
simply connected manifold $M$, and let  $L_f$ be defined by the
associated pullback diagram,
$$
 \xymatrix{
 L_f\dto\rto
& LM\dto \\
 Z\rto^{f}& M\, .
 }
$$
 The shifted homology, $
\bH_\bullet(L_f):=H_{\bullet + p}(L_f)$, is naturally
 an associative  (but, in general, non-commutative)
algebra with respect to the obvious analogue, $\Cap$, of the
Chas-Sullivan product.

\sip

 If $Z$ is a
point, then $L_f\simeq \OM$ and $(\bH_\bullet(L_f),\Cap)=
(H_\bullet(\OM),\star)$.  Chen \cite{C} has
found a nice model for the algebra $(H_\bullet(\OM),\star)$ as the
homology of the free differential algebra $(\R[X], \eth)$, where
$R[X]:=\ot^\bullet \tilde{H}_\bullet(M)[1]$ is the tensor
algebra generated by the reduced  homology (with shifted degree) of $M$  and
$\eth$ is a  differential on $R[X]$ which can be computed by a certain
iterative procedure (see Lemma~2.1 below). In this paper we
generalize Chen's model of $(H_\bullet(\OM),\star)$ to the case of the
Chas-Sullivan  algebra, $(\bH_\bullet(L_f),\Cap)$, for an arbitrary
smooth map $f:Z\rar M$. More precisely, we establish two
isomorphisms of algebras,
$$
\left( \bH_\bullet(L_f),  \Cap\right) =\left(
\Lambda_Z\ot \R[X], d_f\right),
$$
and
$$
\left(
\Lambda_Z\ot \R[X], d_f\right) =
\left(
{\mathrm Hoch}^\bullet\left(\LM,\Lambda_Z\right),
{\mathrm Hochschild\ product}\right),
$$
where $\Lambda_Z$ is the De Rham algebra of $Z$ and $d_f$ is a certain
twisting of the natural differential $d+\eth$ on $\Lambda_Z\ot \R[X]$.
In the case $f=\Id: M\rar M$ the above two results imply $(**)$.
In the case $f: {\mathrm point}\rar M$ the same results reproduce
the Chen model for the homology algebra of the based loop space.

\bip

The paper is organized as follows: In Sections 1 and 2 we review
basic facts of Chen's theory of iterated integrals and of a formal
power series connection. In Section 3 we use the latter to give a
suitable model for the Hochschild cohomology of a dg algebra with
finite dimensional cohomology and then illustrate its work in two
standard examples.
In Section 4 we use the theory of iterated integrals to give a new
proof   of the isomorphism theorem $(*)$. The sections 6-10 are
devoted to a new proof of the isomorphism of algebras theorem
$(**)$. In the final section 11 we discuss the de Rham model for
the string topology on a ``brane" $f:Z\rar M$.

\bip

\sip
{\bf 1. Iterated integrals.} Let $M$ be  a smooth $n$-dimensional
manifold and $PM$ the space of piecewise smooth paths $\ga:
[0,1]\rar M$. For each  simplex, $k=0,1,2,\ldots$,
$$
\Delta^k:= \left\{(t_1,\ldots,t_k)\in \R^k\mid 0\leq t_1\leq t_2
\leq \ldots \leq t_k\leq 1\right\},
$$
consider a diagram,
$$
 \xymatrix{
 PM\times \Delta^k \ar[r]^{\ \ ev_k} \ar[d]_{p} & M^{k+2}  \\
 PM  &
 }
 $$
where $p$ is just the projection on the first factor and $ev_k$ is
the evaluation map,
$$
 \begin{array}{rccc}
   ev_k: &  PM\times \Delta^k  & \lon &   M^{k+2}\vspace{2mm}\\
    & \ga\times (t_1,\ldots,t_k) & \lon  &
   (\ga(0), \ga(t_1),\ldots,\ga(t_k),\ga(1)) .
    \\
 \end{array}
 $$

According to Chen \cite{C}, the {\em space of iterated integrals},
$Ch(PM)$, on $PM$ is a subspace of the de Rham space,
$\Lambda_{PM}$, of smooth differential forms on $PM$. By
definition, $Ch(PM)$ is the image of the following composition of
pull-back and push-forward maps,
$$
\bigoplus_{k=0}^\infty p_*\circ ev_k^*: \ \ \
\bigoplus_{k=0}^\infty \LM\ot \LM^{\ot k}\ot \LM \lon
\Lambda_{PM}.
$$
Chen reserved a special symbol,
$$
\int w_1w_2\ldots w_k :=p_*\circ ev_k^*(1\ot w_1\ot
w_2\ot\ldots\ot  w_k\ot 1),
$$
 for the restriction of $p_*\circ ev_k^*$ to the subspace
  $1\ot \LM^{\ot k}\ot 1 \subset \LM\ot \LM^{\ot k}\ot \LM $.
  The latter is a differential form on $PM$ of degree $|w_1|+ \ldots + |w_k| - k$
  where $|w_i|$ stands for the degree of the form $w_i\in \LM$, $i=1,\ldots,k$.
  If we denote the composition of $ev_k$ with the projection, $M^{k+2}\rar M$,
  to the first (respectively, last) factor by $\pi_0$ (respectively, $\pi_1$), then we
  can write,
  $$
Ch(PM)={\mathrm span}\left\{\pi_0^*w_0\wedge \int w_1w_2\ldots w_k
\wedge \pi_1^*w_{k+1} \right\},
  $$
where $\{w_i\}_{0\leq i \leq k+1}$ are arbitrary differential
forms on $M$.

\bip

The Stoke's theorem,
$$
dp_* = p_*d +  p_*\mid_{\p \Delta^k},
$$
implies \cite{C} \Beqrn d\int w_1w_2\ldots w_k &=& \sum_{i=1}^k
(-1)^{|w_1|+ \ldots + |w_{i-1}| - i}
\int w_1\ldots dw_i\ldots w_k \\
&& -\, \sum_{i=1}^{k-1} (-1)^{|w_1|+ \ldots + |w_{i}| - i}
\int w_1\ldots (w_i\wedge w_{i+1})\ldots w_k \\
&&
-\, \pi_0^*w_1\wedge \int w_2\ldots w_k \\
&& + (-1)^{|w_1|+ \ldots + |w_{k-1}| - k+1}\int w_1\ldots
w_{k-1}\wedge\pi_1^*w_{k}. \Eeqrn The wedge product of iterated
integrals is again an iterated integral \cite{C},
$$
\int w_1\ldots w_k \wedge \int w_{k+1}\ldots
w_{k+l}=\sum_{\sigma\in Sh(k,l)} (-1)^\sigma\int w_{\sigma(1)}
w_{\sigma(2)} \ldots w_{\sigma(k+l)},
$$
where the summation goes over the set, $ Sh(k,l)$, of all shuffle
permutations of type $(k,l)$, and $(-1)^\sigma$ is the Koszul sign
computed as usually but with the assumption that the symbols $w_i$
have shifted degree $|w_i|+1$.

\bip

This all means that $Ch(PM)$ is a dg subalgebra of the de Rham
algebra $\Lambda_{PM}$. In fact, this subalgebra (together with
variants, $Ch(LM):= Ch(PM)\mid_{LM}$ and $Ch(\OM):=
Ch(PM)\mid_{\OM}$) encodes much of the essential information about
the free path space $PM$ and its more important subspaces, $LM$
and $\OM$, of free loops and based loops respectively:

\bip

{\bf 1.1. Fact} \cite{C}. {\em If $M$ is simply connected, then
the cohomology of $Ch(\OM)$ (respectively, $Ch(LM)$) equals the de
Rham cohomology $H^\bullet(\OM)$ (respectively, $H^\bullet(LM)$).}

\bip

\sip

{\bf 2. Formal power series connection and holonomy map.} Let
$(A,d)$ be a unital dg algebra over a field $k$ with finite
dimensional cohomology, say $\dim H^\bullet(A)=n+1$. Let
$\{[1],[e^i]\}_{1\leq i \leq n}$ be a homogeneous basis of the
graded vector space $H^\bullet(A)$, and $\{1^*, x_i\}_{1\leq i
\leq n}$ the associated dual basis of the shifted dual vector
space $\Hom_k(H^\bullet(A),k)[1]$. Let $k\langle
X\rangle:=k\langle x_1,\ldots,x_n \rangle$ be the free graded
associative algebra generated by non-commutative\footnote{As
opposite to $k\langle x_1,\ldots,x_n \rangle$, the free graded
{\em commutative}\, associative algebra generated by homogeneous
indereterminants $x_1,\ldots,x_n$ will be denoted by\
 $k[x_1,\ldots,x_n]$.}
  indeterminates  $x_1,\ldots,x_n$, and
 $k\langle\langle X\rangle\rangle:=k\langle\langle x_1,\ldots,x_n \rangle\rangle$
its formal completion. Finally, let
 $\{e^i\}$ be arbitrary lifts of the cohomology classes $\{[e^i]\}$ to cycles in $A$.

\bip

{\bf 2.1. Lemma} \cite{C}. {\em There exists a degree one element
$\omega$ in the algebra $A\ot k\langle\langle X\rangle\rangle$ and
a degree one differential, $\eth$, of the algebra $k\langle\langle
X\rangle\rangle$ such that \Bi
\item $\eth I\subset I^2$, where $I$ is the maximal ideal in $k\langle\langle X\rangle\rangle$;
\item $\omega \bmod A\ot I^2 = \sum_{i=1}^n e^i\ot x_i$;
\item the Maurer-Cartan equation,
$$
d\omega + \eth\om + \omega\om =0,
$$
is satisfied. \Ei } The proof goes by induction in the tensor
powers of $I$ (cf.\ \cite{C,H}).

\bip

{\bf 2.2. Remark}. If $(A,d)$ is a non-negatively graded connected
and simply connected (i.e.\ $H^0(A,d)=k$ and $H^1(A,d)=0$)
differential algebra, then the above result holds true with
$\om\in A\ot k\langle X\rangle\subset A\ot k\langle\langle
X\rangle\rangle$. This  case is of major interest to us
 as it covers the particular example when $A$ is the De Rham algebra $\LM$
of a connected and simply connected compact manifold. A geometric
significance of the differential $\eth$ for $A=\LM$ was discovered
by Chen in the following beautiful

\bip

{\bf 2.3. Fact} \cite{C}. {\em Let $M$ be a compact simply
connected manifold. Then there is a canonical isomorphism of
algebras,
$$
\left(H_\bullet(\OM), {\mathrm Pontrjagin\ product}\ \star\
\right) =H^\bullet(\R\langle X\rangle,\eth).
$$
}

\bip

Chen's  proof of this statement employs the notion of the {\em
holonomy map},
$$
 \begin{array}{rccc}
   Hol: &  (C_\bullet(\OM),\p) & \lon &   (\R\langle X\rangle,\eth) \\
    & {\mathrm simplex}\ \al: \Delta^k\rar \OM  & \lon  &
   (-1)^{|\al|}\int_{\Delta^k} \al^*(T) ,
    \\
 \end{array}
 $$
from the complex of (negatively graded) smooth singular  chains in
the based loop space $\OM$ to the complex $(\R\langle
X\rangle,\eth)$. The centerpiece  is the restriction to $\OM$ of
the degree zero element, $T\in \Lambda_{PM}\ot k\langle\langle
x_1,\ldots,x_n\rangle\rangle$, given by the following iterated
integral
$$
T:= 1 + \int\om + \int \om\om + \int \om\om\om  + \ldots .
$$
This element is called the {\em transport}\, of a formal power
series connection $\om$ from Lemma~2.1, and has two useful
properties \cite{C}: \Bi
\item $dT + \eth T + \pi_0^*(\om)\wedge T - T\wedge \pi_1^*(\om) =0$, and
\item given any two smooth simplices, $\al_1: \Delta^k\rar PM$
and $\al_2: \Delta^l\rar PM$, such that for any $v\in \Delta^k,
w\in \Delta^l$ the path $\al_2(w)$ begins where the path
$\al_1(v)$ ends (so that the Pontrjagin product $\al_1\star \al_2:
\Delta^k\times \Delta^l\rar PM$ makes sense), then
$$
\int_{\Delta^k\times \Delta^l} (\al_1\times \al_2)^*(T) =
\int_{\Delta^k} \al_1^*(T)
 \int_{\Delta^l} \al_2^*(T).
$$
\Ei

As $T$ restricted to $\OM$ satisfies $dT=-\eth T$, it is now
obvious that $ Hol:   (C_*(\OM),\p)  \lon   (\R\langle
X\rangle,\eth)$ is a morphism of dg algebras. Some extra work
involving Adam's cobar construction shows
 that $Hol$ is a quasi-isomorphism; see \cite{C} for details.

\bip

{\bf 2.4. Remark.} The dual of the finitely generated
 dg algebra $(k\langle\langle X\rangle\rangle,\eth)$
is a dg coalgebra, i.e.\ an $A_\infty$-structure on the vector
space $H^\bullet(A)$. The latter is precisely a
cominimal\footnote{An $A_\infty$ structure on a vector space $V$,
that is a codifferential $Q$ of the free coalgebra ${\ot^\bullet}
V[1]$, is called {\em cominimal} \, if the restriction of $Q$ to
the linear bit $V[1]\subset \ot^\bullet V[1]$ vanishes. Any
$A_\infty$ algebra over a field of characteristic zero can be
represented as a direct sum of a cominimal and a contractible
ones, the former being determined uniquely up to an isomorphism.}
 model for the original dg algebra $(A,d)$ and Chen's formal power series
connection $\om\in A\ot k\langle\langle X\rangle\rangle$ is
nothing but an $A_\infty$ morphism
$$
(H^\bullet(A), {\mathrm cominimal}\ A_\infty\ {\mathrm
structure})\lon (A,d).
$$
Thus Chen's Lemma~2.1 is tantamount to saying that the cohomology,
if finite dimensional,
 of a dg $k$-algebra
has an induced structure of an $A_\infty$-algebra. A fact proved
later in greater generality and independently by Kadeishvili
\cite{K}.

\bip

{\bf 2.5. Remark.} A formal power series connection $\om$ as
described in Lemma~2.1 is not canonical but depends on a number of
choices. There  is, however,  a canonical object underlying this
power series which is best described using the language of
non-commutative differential geometry.

\sip

First we replace the pair $(k\langle\langle
x_1,\ldots,x_n\rangle\rangle,\eth)$ by a germ of the
$n$-dimensional formal non-commutative smooth graded manifold $X$
equipped with a degree one smooth vector field $\eth$. This means
essentially  that we enlarge the automorphism group of
$(k\langle\langle x_1,\ldots,x_n\rangle\rangle,\eth)$ from
$GL(n,k)$ (a change of basis $x_i\rar \sum a_{ij}x_j$ in
$H(A)[1]^*$) to arbitrary formal diffeomorphisms ($x_i\rar \sum
a_{ij}x_j + a_{ijk}x_jx_k+\ldots$), i.e.\ we forget the flat
structure in $H(A)[1]^*$.

\sip

Next we observe that $\om$ defines a flat connection,
$\nabla_\eth\simeq \eth + d + \om$, along the vector field $\eth$
in the trivial bundle over $X$ with typical fibre $A$. It is this
flat $\eth$-connection $\nabla_\om$ which is defined canonically
and invariantly. Chen's formal power series $\om$ is nothing but a
representative of $\nabla_\om$ in a particular coordinate chart
$(x_1,\ldots,x_n)$ on $X$.

\sip

This whole paper (more precisely, everything starting with
Theorem~3.1 below) should have been written using the geometric
language of flat $\eth$-connections\footnote{Let $(X,\eth)$ be a smooth dg manifold,
that is, a graded
manifold $X$ equipped with a degree one smooth vector field $\eth$ such that
$[\eth,\eth]=0$. Let $E$ be a vector bundle over $X$ which we understand
as  a locally free sheaf of $\f_X$-modules, $\f_X$ being the structure sheaf on $X$. An
$\eth$-{\em connection}\, on $E$ is a linear map $\nabla_\eth:E\rar E$ such that
$\nabla_\eth(fe)= (-1)^{|f|}f\nabla_\eth(e) + (\eth f)e$ for any $f\in \f_X$, $e\in E$.
 It is called {\em flat}\,
if $\nabla_\eth^2=0$.       } $\nabla_\eth$   on the non-commutative dg
manifold $(X,\eth)$. After some hesitation we decided not to do it
and work throughout in a particular coordinate patch
$(x_1,\ldots,x_n)$; this makes our formulae more transparent but
at a price of using a non-canonical, a choice of coordinates
dependent
 object $\om$. Of course, ultimately nothing depends on that choice.

\sip

\bip
{\bf  3. A model for Hochschild cohomology.} Let $(A,d)$ be a
unital dg $k$-algebra with finite dimensional cohomology.
Lemma~2.1 says that there is a canonical (see the remark just
above)
 deformation, $ d + \eth + \om$, of the differential
$d$ in $A\ot k\langle\langle X\rangle\rangle$. We shall need below
its Lie bracket version,
$$
d_\om a := da + \eth a + [\om, a], \ \ \ \ \forall a\in A\ot
k\langle\langle X\rangle\rangle.
$$
Clearly, the Maurer-Cartan equation implies $d_\om^2=0$. If
$(M,d)$ is a dg bi-module over $(A,d)$, then formally the same
expression as above defines a differential $d^M_\om$ in $M\ot
k\langle\langle X\rangle\rangle$.

\sip

The following result should be well-known though the author was
not able to trace a reference.

\bip

{\bf 3.1. Theorem.} (i) {\em The  Hochschild cohomology
$Hoch^\bullet(A,A)$ is canonically isomorphic as an associative
algebra to the cohomology $H^\bullet(A\ot k\langle\langle
X\rangle\rangle, d_\om)$.}

\sip

\hspace{27mm} (ii) {\em The  Hochschild cohomology
$Hoch^\bullet(A,M)$ is canonically isomorphic as a vector space to
the cohomology $H^\bullet(M\ot k\langle\langle X\rangle\rangle,
d^M_\om)$.}

\bip

\Proof  We show the proof of (i) only (the proof of statement (ii)
is analogous). The idea is simple: we shall construct a continuous
morphism of topological
  dg algebras,
$$
\left( C^\bullet(A,A):= \Pi_{n\geq 1}Hom_k(A^{\ot n},A)[-n],
\cup_{Hoch},d_{Hoch}
 \right) \lon \left (A\ot k\langle\langle X\rangle\rangle,
{product\ot product},  d_\om\right),
$$
which induces an isomorphism in cohomology.

\bip

Usually $ C^\bullet(A,A)$ gets identified with the vector space of
coderivations of the bar construction on $A$ as it nicely exposes
the dg Lie algebra structure. We are interested, however, in the
dg associative algebra structure on  $ C^\bullet(A,A)$ and thus
have to look for something else.

\sip

 Let $A_1$ and $A_2$ be $A_\infty$-algebras, that is, a pair of co-free
codifferential coalgebras $(B(A_1),\Delta,Q_1)$ and
$(B(A_2),\Delta,Q_2)$, where $B$ stands for the bar construction,
$$
B(A_i):= \bigoplus_{n\geq 1} (A_i[1])^{\ot n},
$$
$\Delta$ for the coproduct,
$$
 \begin{array}{rccl}
   \Delta: &  B(A_i) & \lon & B(A_i)\ot B(A_i)
   \vspace{2mm}\\
    & [a_1|a_2|\ldots|a_n] & \lon  & \sum_{k=1}^{n-1} [a_1,\ldots,a_k]\ot
  [a_{k+1}|\ldots|a_n],
 \end{array}
 $$
and $Q_i: B(A_i)\rar B(A_i)$ for a degree 1 coderivation of
$(B(A_i),\Delta)$ satisfying the equation $Q_i^2=0$. For example,
if, say, $A_2$ is just a dg algebra (which we now assume from now
on) with the differential $d:A_2\rar A_2$ and the product
$\mu_2:A_2\ot A_2\rar A_2$, then the associated codifferential is
given by \Beqrn Q_2[a_1|a_2|\ldots|a_n]&:=&\sum_{k=1}^n
(-1)^{|a_1|+ \ldots + |a_{k-1}| - k}
[a_1|\ldots|da_k|\ldots|a_n] \\
&& -\, \sum_{k=1}^{n-1} (-1)^{|a_1|+ \ldots + |a_{k}| - k}
[a_1|\ldots|\mu_2(a_k,a_{k+1})|\ldots|a_n]. \Eeqrn

\bip

The  vector space of all  coalgebra maps,
$$
Comap(A_1,A_2):={\Hom}_{Coalg}(B(A_1),B(A_2)),
$$
is naturally isomorphic to ${\Hom}_k(B(A_1), A_2[1])=
 \Pi_{n\geq 1}Hom_k(A_1^{\ot n},A_2)[1-n]$.
Then one can turn $Comap(A_1,A_2)[-1]$ into an associative algebra
with the product given by (see e.g.\ \cite{HMS})
$$
f\cup g:= \mu_2\circ (f\ot g)\circ \Delta, \ \ \ \forall f,g\in
{\Hom}_k(B(A_1), A_2).
$$

\bip

A remarkable fact is that every  $A_\infty$ morphism, i.e.\ every
element of the set
$$
\Hom_{A_\infty}(A_1,A_2):=\left\{\om\in{\Hom}_{Coalg}(B(A_1),B(A_2)):
|\om|=0, \om\circ Q_1=Q_2\circ \om\right\},
$$
makes $(Comap(A_1,A_2)[-1],\cup)$ into a dg algebra. Indeed, the
defining equation $\om\circ Q_1=Q_2\circ \om$ for $\om$
(understood  as a degree 1 element of $\Hom_k(B(A_1), A_2)$) takes
the form of a twisting cochain equation,
$$
d\circ \om + \om\circ Q_1 + \om\cup\om\equiv D\om +\om\cup\om=0,
$$
where $D$ is the natural differential in
$\Hom_{complexes}((B(A_1),Q_1), (A_2,d))$.

\bip

Then $ \left(Comap(A_1,A_2)[-1], \cup, \delta_\om\right)$ with the
differential defined by
$$
\delta_\om (f):=Df+ \om\cup f - (-1)^{|f|} f\cup \om
$$
is obviously a topological  dg algebra.

\bip

For example, for any dg algebra $A$, $(Comap(A,A)[-1], \cup,
\delta_{Id})$ is nothing but the
 Hochschild dg algebra  of  $A$.

\bip

Now we assume, in the above notations,
 that $A_2$ is the dg algebra of Theorem~3.1 and $A_1$ is its cohomology
$H(A,d)$ equipped  with an induced $A_\infty$-structure $Q$. As
$H(A,d)$ is finite-dimensional, the dualization of the standard
completion of $(B(H(A,d)),Q)$ is a free topological
 dg algebra
which is nothing but $(k\langle\langle X\rangle\rangle,\eth)$.
Thus the dg algebra,
$$
(A\ot k\langle\langle X\rangle\rangle, {product\ot product},
d_\om),
$$
is canonically isomorphic to $(Comap(H(A,d),A)[-1], \cup,
\delta_\om)$. On the other hand, the composition with the
$A_\infty$-morphism $\om\in Comap(H(A,d),A)$ induces a natural
continuous  map
$$
\om^*: Comap(A,A)[-1]\lon Comap(H(A,d),A)[-1].
$$
A straightforward calculation shows that $\om^*$ is actually a
continuous
 morphism
of topological  dg algebras,
$$
(Comap(A,A)[-1], \cup, \delta_{Id}) \lon (Comap(H(A,d),A)[-1],
\cup, \delta_\om),
$$
that is, a continuous morphism  of the topological dg algebras
$$
\left(C^\bullet(A,A), \cup_{\mathrm Hoch}, d_{\mathrm Hoch}
\right)
\stackrel{\om^*}{\lon } (A\ot k\langle\langle X\rangle\rangle,
{product\ot product}, d_\om).
$$
Being continuous, this morphism preserves the decreasing
filtrations, $C^{\geq r}(A,A)$ and $A\ot I^r$, $I$ being the
maximal ideal in $ k\langle\langle X\rangle\rangle$, and, as it is
easy to see, induces an algebra isomorphism of the $E_1$ terms of
the associated spectral sequences. As  filtrations are complete
and exhaustive, the
 spectral sequences converge to $Hoch^\bullet(A,A)$ and
$H^\bullet(A\ot k\langle\langle X\rangle\rangle, d_\om)$
respectively. Then the Comparison Theorem establishes the required
isomorphism. \hfill $\Box$

\bip

{\bf 3.2. Remark.} If $A$ is connected and simply connected, then
 $k\langle\langle X\rangle\rangle$
can be replaced in the above theorem by $k\langle X \rangle$.

\bip

{\bf 3.3. Example.} Let $M$ be an $n$-dimensional sphere $S^n$.
Let $A=\Bbb R[\nu]$ be a dg subalgebra of the de Rham algebra
$\LM$  generated over $\R$ by a volume form $\nu$. As $A$ is
quasi-isomorphic to $\LM$, $Hoch^\bullet(A,A)=
Hoch^\bullet(\LM,\LM)$.

\sip

Clearly, the data $\eth=0$, $\om=\nu\ot x$, $|x|=1-n$, is a
solution of the Maurer-Cartan equation for $A$,
$$
d\omega + \eth\om + \omega\om =0.
$$
 Hence, by Theorem~3.1, one gets immediately
\Beqrn Hoch^\bullet(\LM,\LM) &=& \frac{\Ker: [\nu\ot x,\ ]:
\R[\nu]\ot\R[x]\lon \R[\nu]\ot\R[x] }
{\Img: [\nu\ot x,\ ]: \R[\nu]\ot\R[x]\lon \R[\nu]\ot\R[x]}\vspace{4mm}\\
 &=&
\left\{ \Ba{ll}
\R[\nu, x]/(\nu^2), \ |\nu|=n, |x|=1-n, & {\mathrm if}\ n \ {\mathrm is\ odd,}\\
\R[\nu,\mu, \tau]/(\nu^2, \mu^2, \nu\mu, \nu\tau ), |\nu|=n,
|\mu|=1, |\tau|=2-2n & {\mathrm if}\ n \ {\mathrm is\ even.} \Ea
\right. \Eeqrn Which is in accordance with the earlier
calculations
 \cite{CJY,FTV}.
The key to the $n$ even case is
$$
\mu\simeq \nu\ot x, \ \ \ \tau\simeq 1\ot x^2.
$$

\bip

{\bf 3.4. Example.} Let $M$ be an $n$-dimensional complex
projective space $\CP^n$. Let $A=\R[h]/h^{n+1}$ be a dg subalgebra
of the de Rham algebra $\LM$  generated over $\R$ by the standard
K\"ahler form $h$.
 As $A$ is quasi-isomorphic to $\LM$, $Hoch^\bullet(A,A)=
Hoch^\bullet(\LM,\LM)$.

\sip

The data,
$$
\om = \sum_{i=1}^n h^{\wedge i}\ot x_i, \ \ \ \ |x_i|=1-2i,
$$
and
$$
\eth=-\sum_{i=2}^n\sum_{j=1}^{i-1}x_jx_{i-j}\frac{\p}{\p x_i},
$$
give a solution of the Maurer-Cartan equations for $A$. Plugging
these into Theorem~3.1, one gets after a  straightforward
 inspection,
$$
Hoch^\bullet(\LM,\LM)= \R[h,\mu,\nu]/(h^{n+1}, h^n\mu,  h^n\nu),
|h|=2, |\mu|=1, |\nu|=-2n.
$$
This is in agreement with a spectral sequence calculation made
earlier in \cite{CJY}. The key to the answer is
$$
h\simeq h\ot 1, \ \ \  \mu\simeq \sum_{i=1}^n ih^{\wedge i}\ot x_i
, \ \ \ \nu\simeq \sum_{i+j=n+1}1\ot x_ix_j.
$$

\bip

\sip

{\bf 4. Holonomy map for the free loop space.} Let $M$ be a
compact smooth manifold, and $\om\in \LM\ot \R\langle\langle
X\rangle\rangle$ a formal power series connection of the de Rham
algebra $(\LM,d)$ (see Lemma~2.1 and Remark~2.4). As $\LM^*$ is
naturally a dg bi-module over $(\LM, d)$, the connection $\om$
gives rise to a differential, $d^*_\om=d^*+\eth + [\om, \ldots]$,
on $\Hom_\R(\LM,\R\langle\langle X\rangle\rangle)$ (cf.\ section
3). Explicitly, for any $f\in \Hom_\R(\LM,\R\langle X\rangle) $
and any $(\ldots)\in \LM$,
$$
(d^*_\om f)(\ldots)=- (-1)^{|f|}f(d\ldots) + \eth f(\ldots) +
\sum_l (-1)^{|f||\om_l|} [p_l\langle X\rangle,
f(w_l\wedge\ldots)],
$$
where we used a decomposition of a formal power series connection
$\om \in \R\langle\langle X\rangle\rangle\ot \LM$ into a sum of
tensor products,
 $\om=\sum_l p_l\langle X\rangle \ot w_l$.

\bip


The transport $T$ of $\om$ when restricted to the free loop
subspace, $LM$, of the path space $PM$ satisfies the equation,
$$
dT + \eth T + \pi^*(\om) T - T\pi^*(\om) =0,
$$
where $\pi: LM\rar M$ is the evaluation map
$\pi_0\mid_{LM}=\pi_1\mid_{LM}$.

\bip

Let $ (C_\bullet(LM),\p)$ be the complex of (negatively graded)
smooth singular  chains in the free loop space.

\bip

{\bf 4.1. Theorem.} {\em The map}
$$
 \begin{array}{rccl}
   Hol: &  \left(C_\bullet(LM),\p\right) & \lon &
   \left(\Hom_\R(\LM,\R\langle\langle X\rangle\rangle), d^*_\om\right) \vspace{2mm}\\
    & {\mathrm simplex}\ \al: \Delta^\bullet\rar LM  & \lon  &
      Hol(\al): \LM \rar \R\langle\langle X\rangle\rangle  \vspace{2mm}\\
     &&&\hspace{13mm} (\ldots) \rar (-1)^{|\al|}
   \int_{\Delta^\bullet} \al^*\left(\pi^*(\ldots)\wedge T\right) ,
    \\
 \end{array}
 $$
{\em is a morphism of complexes. Moreover, it induces an
isomorphism in cohomology if $M$ is simply connected}.

\sip

\Proof That $Hol$ commutes with the differentials is nearly
obvious, \Beqrn Hol(\p \al)&=&
(-1)^{|\al|-1}\int_{\p\Delta^\bullet} \al^*\left(\pi^*(\ldots)
\wedge T\right)\\
&=& (-1)^{|\al|-1}\int_{\Delta^\bullet} d\al^*\left(\pi^*(\ldots)\wedge T\right)\\
&=& (-1)^{|\al|}\int_{\Delta^\bullet}
\al^*\left(\pi^*(-d\ldots)\wedge T
+ \pi^*(\ldots)\wedge (\eth T + \pi^*(\om) T - T\pi^*(\om)\right)\\
&=& d^*_\om Hol(\al). \Eeqrn

For the second part, consider two increasing filtrations,
$$
\cF_p\Hom_\R(\LM,\R\langle\langle X\rangle\rangle):= \left\{ f\in
\Hom_\R(\LM,\R\langle\langle X\rangle\rangle) \ \mid f(w)=0\
\forall w\in \Lambda^{\geq p}M\right\}
$$
and
$$
\cF_p C_*(LM):= \left\{ \Ba{l} \mbox{simplexes $\Delta^\bullet\rar
LM$ whose
composition with $\pi$}\\
\mbox{ land in the skeleton $M^p$}\\
\mbox{for an appropriate smooth triangulation of $M$}\\
 \Ea
 \right\}
$$
Clearly, the map $Hol$ preserves the filtrations. The Leray-Serre
spectral sequence associated with the filtration $\cF_p C_*(LM)$
has
$$
E^2_{\bullet,\bullet}=H_\bullet(M)\ot H_\bullet(\OM).
$$

As $M$ is simply connected, the formal power series connection $w=
\sum_l p_l(x_1,\ldots,x_n)\ot w_l$ has all $w_l\in \LM^{\geq 2}$.
Then using Fact~2.3, it is easy to compute the $E^2$ term of the
spectral sequence associated to $\cF_p\Hom_\R(\LM,\R\langle\langle
X\rangle\rangle)$,
$$
\hat{E}^2_{\bullet,\bullet}=(H^\bullet(M))^*\ot H_\bullet(\OM).
$$

Thus $Hol$ provides an  isomorphism of the $E^2$ terms. Then the
standard spectral sequences comparison argument finishes the
proof. \hfill $\Box$

\bip

The above Theorem in conjunction with Theorem~3.1(ii) immediately
implies the  Jones isomorphism  \cite{J}:

\bip

{\bf 4.2. Corollary.} {\em If $M$ is simply connected, the
holonomy map induces the isomorphism},
$$
H_\bullet(LM, \R) = {\mathrm
Ext}_{\LM\ot\LM}^\bullet\left(\LM,\LM^*\right).
$$

\bip

\sip

{\bf 5. BV operator and an odd holonomy map.} The holonomy map
$Hol$ can also
 be understood as a degree 0 morphism of complexes,
$$
 \begin{array}{rccl}
   \widetilde{Hol}: &  (\LM, d) & \lon &
\left( \Lambda_{LM}\ot \R\langle\langle X\rangle\rangle,
d_{\pi^*(\om)}:= d + \eth + [\pi^*(\om),\ldots]\right)
    \vspace{2mm}\\
    &  \Xi & \lon  & \pi^*(\Xi)\wedge T
    \\
 \end{array}
 $$
\sip

Consider now a degree -1 linear map
$$
 \begin{array}{rccl}
   \widetilde{\Psi}: &  \LM\ot\R\langle\langle  X\rangle\rangle
   & \lon & \Lambda_{LM}\ot\R\langle\langle X\rangle\rangle
    \vspace{2mm}\\
    &  \Xi & \lon  & \sum_{n=0}^\infty\sum_{i=0}^n \int \underbrace{\om\ldots\om}_i
    \Xi \underbrace{\om\ldots\om}_{n-i}
    \\
 \end{array}
 $$

\bip
{\bf 5.1. Lemma.} {\em For any $\Xi\in \LM\ot\R\langle\langle
X\rangle\rangle$},
$$
[T,\pi^*(\Xi)]=d_{\pi^*(\om)} \widetilde{\Psi}(\Xi) +
\widetilde{\Psi}(d_\om\Xi).
$$

Proof is a straightforward calculation of the r.h.s.

\bip

Note that $(\LM\ot 1, d)$ is naturally a subcomplex of
$(\LM\ot\R\langle\langle X\rangle\rangle, d_\om)$.

\bip

{\bf 5.2. Corollary.} {\em The map $ \widetilde{\Psi}_0:=
\widetilde{\Psi}|_{\LM\ot 1}$ induces a degree $-1$ morphism of
complexes,
$$
\left(\LM, d\right) \lon \left(\Lambda_{LM}\ot\R\langle\langle
X\rangle\rangle, d_{\pi^*(\om)}\right),
$$
and, through pairing with chains,  a degree $-1$ morphism of
complexes,}
$$
 \begin{array}{rccl}
   \Psi: &  \left(C_\bullet(LM),\p\right) & \lon &
   \left(\Hom_\R(\LM,\R\langle\langle X\rangle\rangle), d^*_\om\right) \vspace{2mm}\\
    & {\mathrm simplex}\ \al: \Delta^\bullet\rar LM  & \lon  &
      \Psi(\al): \LM \rar \R\langle\langle X\rangle\rangle  \vspace{2mm}\\
     &&&\hspace{11mm} (\ldots) \rar
   \int_{\Delta^\bullet} \al^*\left( \widetilde{\Psi}_0(\ldots)\right) .
    \\
 \end{array}
 $$

\bip

There is a natural action of the circle group $S^1$ on the free
loop space,
$$
Rot: LM\times S^1 \rar LM,
$$
 by rotating the loops, $Rot: (\ga(t),s)\rar \ga(t+s)$.
This action induces a degree $-1$ operator on (negatively graded)
chains,

$$
 \begin{array}{rccl}
   \Delta_{BV} : &  C_\bullet(LM) & \lon &  C_{\bullet+1}(LM)
   \vspace{2mm}\\
    &  \al : \Delta^\bullet\rar LM  & \lon  &  \Delta_{BV}(\al):
      \Delta^\bullet \times S^1 \lon LM  \vspace{2mm}\\
     &&&\hspace{13mm} (z,s) \rar Rot(\al(z),s) .
    \\
 \end{array}
 $$

This operator commutes with the differential $\p$ and hence
induces an operator on singular homology, $\Delta_{BV}:
H_\bullet(LM) \rar  H_{\bullet+1}(LM)$,
 which
satisfies the condition $\Delta^2=0$ as its  iteration
$\Delta_{BV}^{k\geq 2}$ increases the geometric dimension of the
input simplex only by one \cite{CS}.

\bip

{\bf 5.3. Proposition.} $Hol\circ \Delta_{BV} = \Psi$.

\bip

\Proof  Let $q:LM\times S^1\rar LM$ be the natural projection. It
was shown in  \cite{GJP} that the $S^1$ rotation affects an
arbitrary   iterated integral as follows,
$$
q_*\circ Rot^*\left( \pi^*(w_0)\wedge \int w_1w_2\ldots w_k\right)
=\sum_{i=0}(-1)^r \int w_i\ldots w_kw_0w_1\ldots w_{i-1},
$$
where $r=(|w_0| + |w_1| +\ldots + |w_{i-1}| -i)( |w_i| +\ldots +
|w_{k}| -k+i$. Applying this formula to $ \pi^*(\ldots)\wedge T $,
$\forall (\ldots)\in \LM$,
 one immediately obtains the required result.
 $\Box$

\bip

\sip

{\bf 6. Poincare duality.} Here is a deformed version of the
classical duality:

\sip

{\bf 6.1. Theorem.} {\em For any smooth compact oriented manifold
$M$ the linear map
$$
 \begin{array}{rccl}
   \cP : &  \left(\LM\ot\R\langle\langle X\rangle\rangle, d_\om\right) & \lon &
   \left(\Hom_\R(\LM,\R\langle\langle X\rangle\rangle), d^*_\om\right)
    \vspace{2mm}\\
&  \Xi & \lon  &
      \cP(\Xi): \LM \rar \langle\langle X\rangle\rangle,
  \vspace{2mm}\\
     &&&\hspace{10mm} (\ldots)\rar
   \int_M \Xi\wedge \left(\ldots\right) ,
    \\
 \end{array}
 $$
is a degree $-n$ quasi-isomorphism of complexes.}

\sip

\Proof Checking commutativity with differentials is  just an
elementary application of the Stokes theorem. The map $\cP$
obviously preserves natural filtrations of both sides by powers of
the maximal ideal in $ \R\langle\langle X\rangle\rangle$. At the
$E_1^{\bullet,\,\bullet}$ level of the associated spectral
sequences this map induces, by the classical Poincare duality, a
degree $-n$ isomorphism,
$$
\cP: H^\bullet(M)\ot\R\langle X\rangle^\bullet \lon
\Hom_\R(H^\bullet(M), \R)\ot \R\langle X\rangle^\bullet.
$$
Hence the spectral sequences comparison theorem delivers the
required result. $\Box$

\bip

{\bf 6.2. Corollary.} If $M$  is simply connected, then, for any
cycle $\al:\Delta^\bullet  \rar LM$ in $ (C_\bullet(LM),\p)$,
there exists a cycle $\Xi_\al\in  \left(\LM\ot\R\langle X\rangle,
d_\om\right)$ such that,
$$
(-1)^{|\al|}
   \int_{\Delta^\bullet} \al^*\left(T\wedge \pi^*(\ldots)\right)
   =\int_M \Xi_\al\wedge \left(\ldots\right) + d_\om H(\ldots),
 $$
for all $(\ldots)\in \LM$ and some $H\in \Hom_\R(\LM,\R\langle
X\rangle)$.

\sip

Our next task is to find an explicit expression of  such a cycle
$\Xi_\al$ in terms of iterated integrals  and a formal power
series connection.

\bip

\sip

{\bf 7. Thom class and a suitable homotopy.} Let $i:M\hook M\times
M$ be the diagonal embedding, and $[U]\in H^n(M\times M, M\times
M\setminus i(M))$ the associated Thom class. For any tubular neighbourhood,
$Tub$, of $i(M)$ in $M\times M$ one can represent $[U]$ by a
closed $n$-form $U\in \Lambda_{M\times M}$ such that \Bi
\item ${\mathrm support}\,U\subset Tub$, and
\item $pr_*(U)=1$,  where $pr:M\times M\rar M$ is the projection
to the first factor and $1$ is the constant function on $M$. \Ei

Let us consider consider a commutative diagram of maps,
$$
 \xymatrix{
 & M\times LM \drto^{\hat{p}_y}\ddrto_{{p}_y}\ddlto_{p_x} & \\
  && LM\dto_\pi   \\
  M&   & M
 }
$$
where $p_x$ and  $\hat{p}_y$ are natural
 projections to, respectively, the first and second factors, and let us
introduce a list of notations, \Beqrn
\widehat{Tub} &:=& (p_x\times p_y)^{-1}(Tub) \subset M\times LM,\\
\hat{T}_y &:=& (\hat{p}_y)^* T \hspace{18mm}\in \Lambda_{LM\times M}\ot R\langle X\rangle,\\
\hat{U}_{x,y} &:=& (p_x\times p_y)^{*}(U)\hspace{7mm}\in \Lambda^n_{M\times LM},\\
\hat{\om}_x &:=& ({p}_x)^* \om \hspace{18mm}\in \Lambda_{LM\times M}\ot R\langle X\rangle,\\
\hat{\om}_y &:=& ({p}_y)^* \om \hspace{18mm}\in \Lambda_{LM\times M}\ot R\langle X\rangle,\\
\Eeqrn (so that a symbol with the hat always stands for the object
living on $M\times LM$). A point in $\widehat{Tub}$ is a pair $(x,
l_y)$ consisting of a point $x\in M$ and a loop $l_y$ in $LM$
based in $y:=\pi(l)$ such that $(x, y)\in Tub$.

\sip

Let $\nabla$ be an arbitrary smooth connection on $M$. For any two
sufficiently close to each other points $x,y\in M$ there is a
unique vector $V\in T_x M$ such that $ y=\exp_x V$. For any $s\in
[0,1]$ define  the path
$$
\Ba{rccc}
P_{s\cdot x,y}: & [0,1] & \lon & M \vspace{2mm}\\
& t &\lon & exp_x((1-s)tV + sV) \Ea
$$
which is just a $[0,1]$-parameterized geodesic from $s\cdot
x:=exp_x(sV)$ to $y=exp_x V$. Such a path allows us in turn to
introduce a  smooth map$$
 \begin{array}{rccl}
   F: &  \widehat{Tub}\times [0,1] & \lon &
   \widehat{Tub} \vspace{2mm}\\
    & (x,l_y)\times s & \lon  & \left(x,(P_{s\cdot x,y} \star l)\star
    P_{s\cdot x,y}^{-1}\right)
 \end{array}
 $$
where $\star$ is the Pontrjagin products of paths.
$$
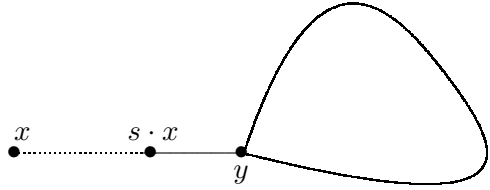
\xy (20,0)*+{};(38,0)*+{} **\crv~pC{~**\dir{.}} ?(0.0)*{\bullet}
*!LD!/^-5pt/{x}
\endxy
\hspace{-2mm} \xy (38,0)*+{};(50,0)*+{} **\crv~pC{}
?(0.0)*{\bullet} *!LD!/^-5pt/{\hspace{-3mm}s\cdot x}
\endxy
\hspace{-2mm}
\xy (50,0)*+{};(50,-1)*+{} **\crv{(90,-10)&(80,10)&(60,30)}
?(0.0)*{\bullet} *!LD!/_-12pt/{y}
\endxy
$$

\bip

\bip

\noindent Note that $F_1=F|_{\widehat{Tub}\times 1}$ is the
identity map $\widehat{Tub}\rar \widehat{Tub}$, while
$F_0=F|_{\widehat{Tub}\times 0}$ factors through the ``based
diagonal",
$$
.
$$
$$
 \begin{array}{ccccccc}
\widehat{Tub}&\hook & M\times LM &\lon & LM &\stackrel{\pi\times
Id}{\hook} & M\times LM
\vspace{2mm}\\
&&(x,l_y) &\lon& (P_{x,y} \star l)\star
    P_{x,y}^{-1}&&
 \end{array}
 $$
\sip

\sip \noindent \sip In particular, the pullback map
$F_0^*:\Lambda_{\widehat{Tub}}\rar \Lambda_{\widehat{Tub}}$
vanishes on any differential form which has a factor of the type
$p_y^*(\ldots) - p_x^*(\ldots)$, $\forall (\ldots)\in \LM$.

\bip

Next we decompose, as usual, the pullback map,
$$
F^* = \ 'F^* \oplus \ ''F^*,
$$
into a bit, $'F^*$,  containing $ds$ and a bit, $'' F^*$, not
containing $ds$ and introduce a homotopy operator,
$$
 \begin{array}{rccl}
   h: &  \Lambda_{\widehat{Tub}}^\bullet & \lon &
    \Lambda_{\widehat{Tub}}^{\bullet-1} \vspace{2mm}\\
    & (\ldots) & \lon  & {\int_0^1}\ {'F^*(\ldots)},
 \end{array}
 $$
such that
$$
\Id = F_0^* + dh + hd.
$$
There are lots of homotopy operators satisfying the above
equation, but the one we
 described above  is particularly suitable:

\bip

{\bf 7.1. Fact.} $ h:  \Lambda_{\widehat{Tub}}^\bullet \rar
    \Lambda_{\widehat{Tub}}^{\bullet-1}$ is a morphism of $p_x^*(\LM)$-modules.

\bip

Extending by zero differential forms with compact support in
$\widehat{Tub}$  to the whole $M\times LM$ we can summarize our
achievements in this section as follows.

\bip

{\bf 7.2. Fact.} There exists a morphism  of $p_x^*(\LM)$-modules,
$$
H=(-1)^n\hat{U}_{x,y}\wedge h: \ \ \  \Lambda_{M\times LM}^\bullet
\lon \Lambda_{M\times LM}^{\bullet+n-1},
$$
such that, for any $(\ldots)\in \LM$ and any
$\Xi\in\Lambda_{M\times LM}$ , \Beqrn \hat{U}_{x,y}\wedge
\Xi\wedge\left( p_y^*(\ldots) - p_x^*(\ldots)\right)
&  = & dH\left(\Xi\wedge\left( p_y^*(\ldots) - p_x^*(\ldots)\right)\right)\\
&& + (-1)^n H\left(d\Xi\wedge\left( p_y^*(\ldots) -
p_x^*(\ldots)\right)\right)
\\
&& + (-1)^{n + |\Xi|} H\left(\Xi\wedge\left( p_y^*(d\ldots) -
p_x^*(d\ldots)\right)\right) \Eeqrn

\bip

In particular, extending $H$ by  $\R\langle\langle
X\rangle\rangle$-linearity to $\Lambda_{M\times
LM}\ot\R\langle\langle X\rangle\rangle$,
 taking $\Xi=\hat{T}_y$ and using the equation, $d\hat{T}_y + \eth\hat{T}_y +
 [\pi_y^*(\om), \hat{T}_y ]=0$, we get
\Beqrn \hat{U}_{x,y}\wedge \hat{T}_y\wedge\left( p_y^*(\ldots) -
p_x^*(\ldots)\right)
&  = & (d+\eth)H\left(\hat{T}_y\wedge\left( p_y^*(\ldots) - p_x^*(\ldots)\right)\right)\\
&& + (-1)^{n-1} H\left([\hat{\om}_y, \hat{T}_y ]\wedge
 \left(p_y^*(\ldots) - p_x^*(\ldots)\right)\right)
\\
&& + (-1)^{n} H\left(\hat{T}_y\wedge\left( p_y^*(d\ldots) -
p_x^*(d\ldots)\right)\right)
\\
&  = &(-1)^{n} H\left([\hat{\om}_y- \hat{\om}_x, \hat{T}_y
]\right)\wedge
 p_x^*(\ldots)\\
&& + (d+\eth)
H\left(\hat{T}_y\wedge\left( p_y^*(\ldots) - p_x^*(\ldots)\right)\right)\\
&& + (-1)^{n-1} H\left([p_y^*(\om\wedge \ldots) - p_y^*(\om\wedge
\ldots), \hat{T}_y ]\right)
\\
&& + (-1)^{n} H\left(\hat{T}_y\wedge\left( p_y^*(d\ldots) -
p_x^*(d\ldots)\right)\right) \Eeqrn

\bip \sip

{\bf 8. An inversion of the deformed Poincare duality map.} We
assume from now on that $M$ is a compact simply connected
manifold. For any cycle $\al:\Delta^\bullet \rar LM$ in $
(C_\bullet(LM),\p)$ there is associated a diagram of maps,
$$
 \xymatrix{
 & M\times \Delta^\bullet\drto^{{\mathrm Id}\times \al}\dlto_{p} & \\
  M&   & M\times LM
 }
$$
 Then,
for any $(\ldots)\in \LM$, we have \Beqrn
 (-1)^{|\al|}  Hol(\al)(\ldots) &=&
   \int_{\Delta^\bullet} \al^*\left(T\wedge \pi^*(\ldots)\right)\\
&=&\int_{M\times \Delta^\bullet} (\Id\times\al)^*
\left(\hat{U}_{x,y}\wedge \hat{T}_y\wedge p_y^*(\ldots)\right)\\
&=& \int_{M\times \Delta^\bullet} (\Id\times\al)^*
\left(\hat{U}_{x,y}\wedge \hat{T}_y\wedge p_x^*(\ldots)\right)\\
&& +\int_{M\times \Delta^\bullet} (\Id\times\al)^*
\left(\hat{U}_{x,y}\wedge \hat{T}_y\wedge(p_y^*(\ldots)-p_x^*(\ldots))\right)\\
\Eeqrn
Hence, setting \Beqrn \Xi_\al &=&  p_*\circ
(\Id\times\al)^*\left(\hat{U}_{x,y}\wedge \hat{T}_y+ (-1)^{n}
H\left([\hat{\om}_y- \hat{\om}_x, \hat{T}_y ]\right)\right)\in
\LM\ot\R\langle X\rangle
\Eeqrn and,
$$
 \begin{array}{rccl}
   H_\al: &  \LM & \lon & \R\langle X\rangle
     \vspace{2mm}\\
    & (\ldots) & \lon  & \int_{M\times \Delta^\bullet}H\left(\hat{T}_y\wedge
    \left( p_y^*(\ldots) - p_x^*(\ldots)\right)\right)
 \end{array}
 $$
we get,
$$
(-1)^{|\al|}Hol(\al)(\ldots) = \int_M \Xi_\al\wedge
\left(\ldots\right) + d^*_\om H_\al(\ldots),
$$
where we have used a calculation at the end of 7.2.


\bip

We almost proved the following

\sip

{\bf 8.1. Theorem.} {\em The degree $n$ linear map
$$
 \begin{array}{ccl}
    (C_\bullet(LM),\p)  & \lon & (\LM\ot\R\langle X\rangle, d_\om)
     \vspace{2mm}\\
     \al & \lon  & (-1)^{|\al|}
     \Xi_\al
 \end{array}
 $$
is a quasi-isomorphism of complexes whose composition with $\cP$
induces on cohomology the map $Hol$}:

$$
 \xymatrix{
   H_\bullet(LM)\ar[dr]_{\hspace{-18mm}Hol}\ar[r]    &
H^{\bullet+n} (\LM\ot\R\langle X\rangle, d_\om))  \ar[d]_{\cP}  \\
&   H^\bullet(Hom_\R(\LM\ot \R\langle X\rangle, d^*_\om)
 }
$$

\bip

\Proof  We need only to check that the map is well defined, i.e.\
that $\Xi_\al$ is a cycle. The latter can be suitably represented
as
$$
\Xi_\al =  p_*\circ (\Id\times\al)^*\left\{\hat{U}_{x,y}\wedge
\left(\hat{T}_y+
 h([\hat{\om}_y- \hat{\om}_x, \hat{T}_y ])\right)\right\}
$$
where we understand $\hat{T}_y+
 h([\hat{\om}_y- \hat{\om}_x, \hat{T}_y ])$ as a  differential form on
 $\widehat{Tub}$ and $\hat{U}_{x,y}\wedge
(\hat{T}_y+
 h([\hat{\om}_y- \hat{\om}_x, \hat{T}_y ]))$ as a differential form
on $M\times LM$ (extended by zero from  $\widehat{Tub}$).

\bip

As opposite to $d_{\hat{\om}_y}\hat{T}_{y}=0$, we have,
\Beqrn d_{\hat{\om}_x}\left\{\hat{T}_y+
 h([\hat{\om}_y- \hat{\om}_x, \hat{T}_y ])\right\}& = &
 -\left[\hat{\om}_y- \hat{\om}_x, \hat{T}_y \right] + (d+ \eth)
 h([\hat{\om}_y- \hat{\om}_x, \hat{T}_y ])
 -\left[\hat{\om}_x,h([\hat{\om}_y- \hat{\om}_x, \hat{T}_y ])\right]
 \\
 &=& -h(d[\hat{\om}_y- \hat{\om}_x, \hat{T}_y ])+
 \eth h([\hat{\om}_y- \hat{\om}_x, \hat{T}_y ])
 + h(\left[\hat{\om}_x,[\hat{\om}_y- \hat{\om}_x, \hat{T}_y ])\right]
 \\
&=&
 h([\hat{\om}_y\hat{\om}_y- \hat{\om}_x\hat{\om}_x, \hat{T}_y ])
 -h([\hat{\om}_y- \hat{\om}_x,[\hat{\om}_y, \hat{T}_y ]])
 -h(\left[\hat{\om}_x,[\hat{\om}_y- \hat{\om}_x, \hat{T}_y ])\right]
 \\
 &=& 0.
\Eeqrn Pushing down this equation along $\al: \Delta^\bullet\rar
LM$ finally gives $d_\om\Xi_\al=0$. \hfill $\Box$

\bip \sip

{\bf 9. Inverse Poincare map revised.} As the calculation in
beginning of Section 8 shows, the class $\Xi_\al =  p_*\circ
(\Id\times\al)^*\left\{\hat{U}_{x,y}\wedge \left(\hat{T}_y+
 h([\hat{\om}_y- \hat{\om}_x, \hat{T}_y ])\right)\right\}$
 provides  us with a natural lift of the holonomy map
$[Hol]: H_\bullet(LM)\rar  H^\bullet(Hom_\R(\LM\ot \R\langle
X\rangle, d^*_\om)$ along  the deformed Poincare duality map,
$[\cP]: H^\bullet(\LM\ot \R\langle X\rangle, d_\om)\rar
 H^\bullet(Hom_\R(\LM\ot \R\langle X\rangle, d^*_\om)$.
However, it is not unique.

\bip

We shall use below another geometrically transparent lift of $
[\cP]^{-1}\circ [Hol]$ from
 the cohomology to the chain level.
Note in this connection that $\Xi_\al$ has to be specified only up
to
 a $d_\om$-exact term implying that the form $\hat{T}_y+
 h([\hat{\om}_y- \hat{\om}_x, \hat{T}_y ])$ on  $\widehat{Tub}$ has to be specified only up to
a $d_{ \hat{\om}_x}$-exact term.

\bip

Recall that for a tubular neighbourhood, $Tub$, of the diagonal in
$M\times M$, the subspace $\widehat{Tub}\subset M\times LM$
consists of pairs $(x, l_y)$, $x$ being a point $x\in M$ and $l_y$
a loop
 in $LM$ based at $y:=\pi(l_y)$, such that $(x,y)\in Tub$. Let
$P_{x,y}$ be the geodesic from $x$ to $y$ for some fixed affine
connection on $M$. Consider a map,
$$
 \begin{array}{rccl}
   i_{x,y}: &  \widehat{Tub} & \lon & LM
    \vspace{2mm}\\
    & (x,l_y) & \lon  & (P_{x,y} \star l_y)\star
    P_{x,y}^{-1}.
 \end{array}
 $$

\bip

{\bf 9.1. Proposition.} $\hat{T}_y+
 h([\hat{\om}_y- \hat{\om}_x, \hat{T}_y])= i_{x,y}^*T \mod \Img d_{ \hat{\om}_x}$.

 \bip

 \Proof Comparing the map $i_{x,y}$ with the homotopy $F$ in Section 7, one notices that
$F_0$ is precisely the composition of $i_{x,y}$ with the natural
embedding,
$$
 \begin{array}{rccl}
   emb:& LM & \lon & M\times LM
    \vspace{2mm}\\
    &l & \lon  & (\pi(l), l).
 \end{array}
 $$
Then applying the equality of maps, $\Id=i_{x,y}\circ emb + dh
+hd$, to $\hat{T}_y$ one obtains, \Beqrn
 i_{x,y}^*T &=& \hat{T}_y - dh(\hat{T}_y))-h(d\hat{T}_y))\\
 &=&  \hat{T}_y + h([\hat{\om}_y- \hat{\om}_x, \hat{T}_y) - d_{ \hat{\om}_x}h(\hat{T}_y).
\Eeqrn \hfill $\Box$

\bip

{\bf 9.2. Corollary.} {\em The degree $n$ linear map
$$
 \begin{array}{rccl}
   \mbox{\bf Hol}: &(C_\bullet(LM),\p)  & \lon & (\LM\ot\R\langle X\rangle, d_\om)
     \vspace{2mm}\\
    & (-1)^{|\al|}\al & \lon  &
      \Ga_\al:= p_*\circ (\Id\times\al)^*\left(\hat{U}_{x,y}\wedge
     i_{x,y}^*T \right)
 \end{array}
 $$
is a quasi-isomorphism of complexes making the following diagram

$$
 \xymatrix{
     &   H^{\bullet+n} (\LM\ot\R\langle X\rangle, d_\om))  \ar[d]_{\cP}  \\
  H_\bullet(LM)\ar[r]_{\hspace{-18mm}Hol}\ar[ur]^{\mathbf Hol}
 &   H^\bullet(Hom_\R(\LM\ot \R\langle X\rangle, d^*_\om)
 }
$$
commutative.}

\bip

\sip

{\bf 10. Compatibility with products.} Now we have all the threads
to show
 a new proof of a remarkable
result of \cite{CJ,Co,Tr} which establishes an isomorphism of {\em
algebras},
$$
\left( \bH_\bullet(LM), \mbox{Chas-Sullivan product}\  \Cap\right)
=\left({\mathrm Ext}_{\LM\ot\LM}^\bullet\left(\LM,\LM\right),
{\mathrm Yoneda\ product}\right).
$$
It follows immediately from Theorem~3.1(i) and the following

\bip

{\bf 10.1. Theorem.} {\em The degree $0$ linear map
$$
 \begin{array}{rccl}
    \mbox{\bf Hol}:& (C_{\bullet+n}(LM),\p)  & \lon & (\LM\ot\R\langle X\rangle, d_\om)
     \vspace{2mm}\\
     &\al & \lon  & (-1)^{|\al|}\Ga_\al
 \end{array}
 $$
induces on cohomology a quasi-isomorphism of algebras},
$$
\left( \bH_\bullet(LM),  \Cap\right) \lon
H^\bullet(\LM\ot\R\langle X\rangle, d_\om).
$$

\bip

\Proof We shall use small Latin letters, ``coordinate variables",
to distinguish  various copies of $M$ and $LM$, e.g.\ $M_x$,
$M_y$, $LM_y$, etc.

\sip

Let $\al_y:\Delta^\bullet\rar LM_y $ and
$\al_z:\Delta_z^\bullet\rar LM_z $ be any two
 cycles in
$(C_{\bullet+n}(LM),\p)$ such that their projections,
$\pi_y\circ\al_y:\Delta_y^\bullet\rar M_y$ and
$\pi_z\circ\al_z:\Delta_z^\bullet\rar M_z$ intersect transversally
at a cycle $\pi_y\circ\al_y\cap \pi_z\circ\al_z:
\Delta_{yz}^\bullet\rar M$, where
$$
\Delta_{yz}^\bullet:=(\pi_y\circ\al_y\times
\pi_z\circ\al_z)^{-1}\left(
\left\{\pi_y\circ\al_y(\Delta_y^\bullet) \times
\pi_z\circ\al_z(\Delta_z^\bullet)\right\}\cap \left \{{\mathrm
diagonal\ in}\ M_y\times M_z \right\}\right).
$$
Chas and Sullivan \cite{CS} define the product chain, $\al_y\Cap
\al_z$, as
$$
\begin{array}{rccl}
   \al_1\Cap \al_2: \Delta_{xy}^\bullet & \lon & LM
     \vspace{2mm}\\
     z & \lon  & \al_1(z)\star \al_2(z).
 \end{array}
$$

The theorem will follow if we show that
$$
\hspace{45mm} \Gamma_{\al_y}\wedge
\Gamma_{\al_z}=\Gamma_{\al_y\Cap \al_z} \mod \Img d_\om.
\hspace{45mm} (\bigstar)
$$
In the obvious notations associated with the commutative diagram,
$$
 \xymatrix{
 & M_x\times \Delta_y^\bullet\times \Delta_z^\bullet
\drto^{{\mathrm Id}\times \al_y\times \al_z}
 \dlto_{p} & \\
  M_x&   & M_x\times LM_y\times LM_z
 }
$$
we can represent the l.h.s. of $(\bigstar)$ as follows
$$
\Gamma_{\al_y}\wedge \Gamma_{\al_z}=
 p_*\circ (\Id_x\times\al_y\times \al_z)^*
 \left(\hat{U}_{x,y}\wedge i^*_{x,y}T\wedge \hat{U}_{x,z}\wedge i^*_{x,z}T\right).
$$

Using homotopy 7.2 we can transform $\ldots\wedge
\hat{U}_{x,z}\wedge\ldots$ into
$\ldots\wedge\hat{U}_{y,z}\wedge\ldots$ modulo $
d_{\hat{\om}_x}$-exact terms: \Beqrn
 \hat{U}_{x,y}\wedge \hat{U}_{x,z}\wedge i^*_{x,y}T\wedge i^*_{x,z}T
&=&\hat{U}_{x,y}\wedge \hat{U}_{y,z}\wedge i^*_{x,y}T\wedge
i^*_{x,z}T
 + d H_{x,y}\left( (\hat{U}_{x,z}- \hat{U}_{y,z})
 \wedge i^*_{x,y}T\wedge i^*_{x,z}T\right)\\
&& + H\left( (\hat{U}_{x,z}- \hat{U}_{y,z})\wedge d(
i^*_{x,y}T\wedge i^*_{x,z}T)\right)
\\
&=& \hat{U}_{x,y}\wedge \hat{U}_{y,z}\wedge \hat{T}_{x,y}\wedge
\hat{T}_{x,z} + d_{\hat{\om}_x}H_{x,y}\left( (\hat{U}_{x,z}-
\hat{U}_{y,z}) \wedge i^*_{x,y}T\wedge i^*_{x,z}T\right). \Eeqrn
Hence
$$
\Ga_{\al_y}\wedge \Ga_{\al_z}=
 p_*\circ (\Id\times\al_y\times \al_z)^*
 \left(\hat{U}_{x,y}\wedge \hat{U}_{y,z}\wedge i^*_{x,y}T\wedge i^*_{x,z}T\right)
\mod \Img d_\om.
$$

Let $Tub_{y,z}$ be a tubular neighbourhood of the diagonal in
$M_y\times M_z$ and let $\widetilde{{Tab}}_{y,z}$ be its pre-image
under the natural evaluation projection,
$$
\pi_y\times \pi_z: LM_y\times LM_z \lon M_y\times M_z.
$$
Similarly, let $Tub_{x,y,z}$ be a tubular neighbourhood
(supporting the form ${U}_{x,y}\wedge {U}_{y,z}$) of the small
diagonal in $M_x\times M_y\times M_z$ and let
$\widetilde{{Tab}}_{x,y,z}$ be its pre-image under the natural
map,
$$
\Id_x\times \pi_y\times \pi_z: M_x\times LM_y\times LM_z \lon
M_x\times  M_y\times M_z.
$$

The form $ \hat{U}_{y,z}\wedge i^*_{x,y}T\wedge i^*_{x,z}T$ is
supported in $M_x\times \widetilde{{Tab}}_{y,z}\subset M_x\times
LM_y\times LM_z $. As $Tub_{y,z}$ is homotopy equivalent to $M_y$,
there exists a one parameter family of maps,
$$
   G:  \widetilde{Tub}_{y,z}\times [0,1] \lon\
   \widetilde{Tub}_{y,z}
 $$
 with the property  that
$G_1=G|_{ \widetilde{Tub}_{y,z}\times 1}$ is the identity map
while $G_0=G|_{\widetilde{Tub}_{y,z}\times 0}$ factors through the
composition
$$
 \widetilde{Tub}_{y,z}  \stackrel{ r}{\lon}
  LM_y\times _{M_y} LM_y
\hook  LM_y\times LM_z
$$
for a smooth map $r$. Then one applies the associated to
$\Id_x\times G$ homotopy equality of linear maps,
$$
\Id = \Id_x\ot (G_0^* + dh_G + h_Gd):
\Lambda_{\widetilde{Tub}_{x,y,z}} \lon
 \Lambda_{\widetilde{Tub}_{x,y,z}},
$$
to the form  $i^*_{x,y}T\wedge i^*_{x,z}T$ and gets
$$
 i^*_{x,y}T\wedge i^*_{x,z}T
=  \hat{i}^*_{x,y} \left( q_1^*T\wedge q_2^*T\right) \mod \Img
d_{\hat{\om}_x},
$$
where $q_1$ and $q_2$ are the natural projections,
$$
 \xymatrix{
 & LM_x\times_{M_x} LM_x\drto^{q_2}\dlto_{q_1} & \\
 LM_x&   & LM_x
 }
$$
and $ \hat{i}_{x,y}$ is given
 by (cf.\ with  $ {i}_{x,y}$ in the beginning of Sect.\ 9)
$$
 \begin{array}{rccl}
   \hat{i}_{x,y}: &  \widehat{Tub}_{x,y,z}\
  & \lon & LM_x\times_{M_x} LM_x
    \vspace{2mm}\\
    & (x,l_y, l_z) & \lon  & (P_{x,y} \star r(l_y,l_z))\star
    P_{x,y}^{-1}.
 \end{array}
 $$

There is a natural Pontrjagin product map
$$
 LM_x\times_{M_x} LM_x \stackrel{\star}{\lon} LM_x,
$$
and the basic property
 of the transport of a formal power
series connection says  that
$$
\star^*T= p_1^*(T)\wedge p_2(T).
$$
Hence we can eventually write
$$
 \hat{U}_{x,y}\wedge \hat{U}_{x,z}\wedge i^*_{x,y}T\wedge i^*_{x,z}T
=  \hat{U}_{x,y}\wedge \hat{U}_{y,z}\wedge \hat{i}^*_{x,y}T \mod
\Img  d_{\hat{\om}_x},
$$
which implies
$$
\Ga_{\al_y}\wedge \Ga_{\al_z}=
 p_*\circ (\Id_x\times\al_y\times \al_z)^*
 \left( \hat{U}_{x,y}\wedge \hat{U}_{y,z}\wedge \hat{i}^*_{x,y}T \right)
\mod \Img d_\om.
$$
Recalling the definition of $\al_y\Cap \al_z$ and using
transversality of $p_y\circ\al_y$ and $p_z\circ\al_z$ to get rid
of the Thom class $U_{y,z}$ through $z$-integration one finally
obtains
$$
\Ga_{\al_y}\wedge \Ga_{\al_z}=
 p_*\circ (\Id_x\times\al_y\Cap \al_z)^*
 \left(\hat{U}_{x,y}\wedge  i^*_{x,y}(T)\right)
\mod \Img d_\om.
$$
Which is precisely the desired result $(\bigstar)$. \hfill $\Box$

\bip

{\bf 11. Loops based at a brane.} If $f:Z\rar M$ is a smooth map from a
compact oriented $p$-dimensional manifold to a compact simply
connected manifold $M$, then one can consider
 a space, $L_f$, defined by the
pullback diagram,
$$
 \xymatrix{
 L_f\dto_\mu\rto^\nu
& LM\dto^\pi \\
 Z\rto^{f}& M\ .
 }
$$
 The shifted homology,  $
\bH_\bullet(L_f):=H_{\bullet + p}(LM)$, is obviously  an
associative algebra with respect to the Chas-Sullivan product
$\Cap$ (see \cite{CS,S}).
If $f$ is an embedding, then $L_f$ is a subspace of $LM$ consisting of loops in $M$ which
are based at $Z \subset M$.

\sip

The
 pull-back map on differential  forms, $f^*: \Lambda_M\rar \Lambda_Z$,
 is naturally extended to the tensor product $\LM\ot\R\langle X\rangle$
 as $f^*\ot\Id$ and is denoted by the same symbol $f^*$. The map
$$
f^*:\left( \LM\ot\R\langle X\rangle, d_\om\right) {\lon}
\left(\Lambda_Z\ot\R\langle X\rangle, d_{f^*(\om)}\right)
$$
is obviously a morphism of dg algebras.

\sip

The proofs of Theorems 4.1, 6.1 and 10.1 apply with trivial
changes to their ``brane"  versions:

\sip

{\bf 11.1. Theorem.} {\em The map}
$$
 \begin{array}{rccl}
   Hol: &  \left(C_\bullet(L_f),\p\right) & \lon &
   \left(\Hom_\R(\Lambda_Z,\R\langle X\rangle), d^*_{f^*(\om)}\right) \vspace{2mm}\\
    & {\mathrm simplex}\ \al: \Delta^\bullet\rar L_f & \lon  &
      Hol(\al): \Lambda_Z \rar \R\langle X\rangle \vspace{2mm}\\
     &&&\hspace{13mm} (\ldots) \rar (-1)^{|\al|}
   \int_{\Delta^\bullet} \al^*\left\{\mu^*(\ldots)\wedge \nu^*(T)\right\} ,
    \\
 \end{array}
 $$
{\em is a quasi-isomorphism of complexes}.

\sip


{\bf 11.2. Theorem.} {\em The linear map
$$
 \begin{array}{rccl}
   \cP : &  \left(\Lambda_Z\ot\R\langle X\rangle, d_{f^*(\om)}\right) & \lon &
   \left(\Hom_\R(\Lambda_Z,\R\langle X\rangle), d^*_{f^*(\om)}\right)
    \vspace{2mm}\\
&  \Xi & \lon  &
      \cP(\Xi): \Lambda_Z \rar \R\langle X\rangle,
  \vspace{2mm}\\
     &&&\hspace{10mm} (\ldots)\rar
   \int_Z \Xi\wedge \left(\ldots\right) ,
    \\
 \end{array}
 $$
is a degree $-p$ quasi-isomorphism of complexes.}

\sip

{\bf 11.3. Theorem.} {\em The degree $p$ linear map
$$
 \begin{array}{rccl}
   \mbox{\bf Hol}: &(C_\bullet(L_f),\p)  & \lon &
(\Lambda_Z\ot\R\langle X\rangle, d_{f^*(\om)})
     \vspace{2mm}\\
    & (-1)^{|\al|}\al & \lon  &
       p_*\circ (f\times\mu\circ\al)^*\left(\hat{U}_{x,y}\wedge
     i_{x,y}^*T \right)
 \end{array}
 $$
associated with the diagram,
$$
 \xymatrix{
 & Z\times \Delta^\bullet\drto^{{f}\times \mu\circ\al}\dlto_{p} & \\
  Z&   & M\times LM\ \ \ ,
 }
$$
is a morphism of complexes inducing on cohomology an isomorphism
of algebras,
$$
\left( \bH_\bullet(L_f),  \Cap\right) \lon
H^\bullet(\Lambda_Z\ot\R\langle X\rangle, d_{f^*(\om)}).
$$
Moreover, if $f$ is an embedding, then the pull-back map $f^*: H^\bullet(\LM\ot\R\langle
X\rangle, d_{\om})\rar  H^\bullet(\Lambda_Z\ot\R\langle X\rangle,
d_{f^*(\om)})$ corresponds to the intersection map $\cap_{[L_f]}:
\bH_\bullet(LM) \rar \bH_\bullet(L_f)$.}

\bip

\sip

{\bf 11.4. Example.} Let  $Z$ be the projective space $\CP^m$
and let $f: Z\rar M=\CP^n$, $m< n$, be a linear embedding. A
straightforward calculation using  Theorem 11.3 gives, \Beqrn
\left(\bH_\bullet(L_f,  \Cap\right) &=&
H^\bullet\left( \R[h]/h^{m+1}\ot \R\langle x_1,\ldots,x_n\rangle,
\
 d_{f^*(\om)}= \eth + \left[\sum_{i=1}^m h^i\ot x_i,\ldots\right]
\right)\\
&=&\R[h,\nu,x]/h^{m+1}, \ |h|=2, |\nu|=-2n, |x|=-1, \Eeqrn where
$\eth$ is the same as in Example 3.4. (The key to the r.h.s. is
$x=\sum_{i=1}^{m+1} i h^{i-1}\ot x_i, \nu\sim \sum_{i+j=n+1}1\ot
x_ix_j$.) Moreover the intersection map
$$
 \xymatrix{
  \bH_\bullet(L\CP^n)\ar@{=}[d] \rto^{\cap_{[L_f]}}
& \bH_\bullet(L_f)\ar@{=}[d] \\
 \R[h,\mu,\nu]/(h^{n+1},h^n\mu,h^n
\nu)\rto& \R[h,x,\nu]/h^{m+1}
 }
$$
is given on generators as follows (see Example 3.4 again for the
notations used in the l.h.s.),
$$
\Ba{ccc}
h &\lon & h\\
\nu &\lon& \nu \\
\mu &\lon & hx . \Ea
$$

\bip

\bip

{\small Acknowledgement. It is a pleasure to thank Ralph Cohen,
Thomas Tradler and especially Dennis Sullivan  for valuable
correspondences and remarks.}

\pagebreak

\bip

\bip

{\small

\bip

\bip

Department of Mathematics

Stockholm University

10691 Stockholm

Sweden


\end{document}